\documentclass[12pt]{amsart}
\usepackage{graphicx}
\newtheorem{thm}{Theorem}[section]
\newtheorem{cor}[thm]{Corollary}
\newtheorem{lem}[thm]{Lemma}

\theoremstyle{definition}
\newtheorem{defn}[thm]{Definition}
\theoremstyle{remark}
\newtheorem{rem}[thm]{Remark}
\theoremstyle{problem}
\newtheorem{prob}[thm]{Problem}
\numberwithin{equation}{section}
\begin{document}

\title[Ulam stability problem]{
On the generalized quadratic mappings in quasi-Banach modules over a $C^*$--algebra
$^\ast$}%

\author[H. Kim]{Hark-Mahn Kim}%
\address{\noindent Department of Mathematics,
 Chungnam National University,
 220 Yuseong-Gu, Daejeon, 305-764,
 Republic of Korea
}%
\email{hmkim@math.cnu.ac.kr}%

\author[D. Lee]{Don O Lee  $^\dag$}%
\address{\noindent Information Center for Mathematical Sciences,
Korea Advanced Institute of Science and Technology, 373-1 Guseong-dong, Yuseong-Gu, Daejeon, 305-701,
 Republic of Korea
}%
\email{dolee@kaist.ac.kr}%

\thanks{
%$^\ast$ This work was supported by the second Brain Korea 21 Project in 2006.
This Article was submitted in The Journal of Mathematical Analysis and Applications \newline
\indent $^\dagger$ Corresponding author:dolee@kaist.ac.kr}%

%\date{}%
%\dedicatory{}%
%\commby{}%
% ----------------------------------------------------------------
\begin{abstract}
Let $n>2$ be a positive integer. In this paper, we obtain the general solution of the following functional equation
\begin{eqnarray*}
 n \sum_{1 \le i<j \le n} Q\left( x_i-x_j \right)=\sum_{i=1}^{n}Q\left(\sum_{j =1}^n x_j -n x_i\right)
\end{eqnarray*}
which is derived from the centroid of the $n$ distinct vectors $x_1, \cdots, x_n$
in an inner product space.
Furthermore, we prove that a mapping $f$ between quasi-Banach modules over a $C^*$-algebra
satisfying approximately the equation
can be approximated by a quadratic mapping $Q$ satisfying exactly the equation such that $\|f(x)-Q(x)\|$ is bounded.
\end{abstract}
\maketitle
% ----------------------------------------------------------------
\section{Introduction}

The stability problem of functional equations originated from a question of S.M. Ulam \cite {ul60}
concerning the stability of group homomorphisms:
``When is it true that by slightly changing the hypotheses of a theorem
one can still assert that the thesis of the theorem remains true or approximately true?''
If the answer is affirmative, then
we would say the equation of homomorphism $H(x \ast y)=H(x)\diamond H(y)$ is stable.
The concept of stability for a functional equation arises when we replace
the functional equation by an inequality which acts as a perturbation of the equation.

First, Ulam's question for approximately additive mappings was solved by D.H. Hyers \cite {hy41}.
In 1951, D.G. Bourgin \cite{bo51} was the second author to treat the Ulam stability problem for additive mappings.
Th.M. Rassias \cite{ra78} succeeded in extending the result of Hyers' theorem by weakening
the condition for the Cauchy difference to be unbounded.
A number of mathematicians were attracted to this result of Th.M. Rassias and stimulated to investigate
the stability problems of functional equations.
%And then, G.L. Forti \cite{fo80} and P. G{\v a}vruta \cite{ga94}
%have generalized the result of Th.M. Rassias' theorem,
%which permitted the Cauchy difference to become arbitrary unbounded.

Now, a square norm on an inner product space satisfies the important parallelogram equality
$%\begin{eqnarray*}
\|x+y\|^2+\|x-y\|^2=2(\|x\|^2 +\|y\|^2 )
$%\end{eqnarray*}
 \ for all vectors $x,y.$
If $\triangle ABC$ is a triangle in a finite dimensional Euclidean space
 and $I$ is the center of the side $\overline{BC},$
then the following identity
$%\begin{eqnarray*}
\|\overrightarrow{AB}\|^2+\|\overrightarrow{AC}\|^2=2(\|\overrightarrow{AI}\|^2+\|\overrightarrow{CI}\|^2)
$%\end{eqnarray*}
\ holds for all vectors $A,B$ and $C.$
The following functional equation, which was motivated by these equations,
\begin{eqnarray}
 Q(x+y)+Q(x-y)=2Q(x)+2Q(y) \label{fe1}
\end{eqnarray}
is called a {\it quadratic functional equation}, and every solution of the
 equation (\ref{fe1}) is said to be a {\it quadratic mapping.}
A Hyers-Ulam stability problem for the quadratic functional equation (\ref{fe1}) was first solved by
F. Skof \cite{sk83}.
C. Borelli and G.L. Forti \cite{bf95} generalized the stability result
of the quadratic functional equation.
The stability problems of several functional equations have been extensively investigated
by a number of authors
and there are many interesting results concerning this problem
 \cite{bjl04,fo04,ga94,hir98,ra00}.
Furthermore, C. Park \cite{pa04} have proved the Hyers-Ulam-Rassias
stability problem for functional equations in Banach modules over a $C^*$-algebra.

%%========================================================================================

Now, if $\triangle XYZ$ is a triangle in a finite dimensional Euclidean space and
$G:=\frac{X+Y+Z}{3}$ is the center of gravity of the triangle,
then a simple direct calculation and the definition of the norm
yields the following identity
\begin{eqnarray}\label{qu-1}
\|\overrightarrow{XY}\|^2+\|\overrightarrow{YZ}\|^2+\|\overrightarrow{ZX}\|^2
=3\left(\|\overrightarrow{XG}\|^2+\|\overrightarrow{YG}\|^2+\|\overrightarrow{ZG}\|^2 \right).
\end{eqnarray}

Employing the above identity (\ref{qu-1}), we introduce the new functional equation,
\begin{eqnarray} \label{fe2}
&& 3Q(x-y)+3Q(y-z)+3Q(x-z)\\\nonumber
&&\qquad =Q(y+z-2x)+Q(x+z-2y)+Q(x+y-2z)
\end{eqnarray}
for a mapping $Q:U \rightarrow V$ and for all vectors $x,y,z \in U,$
where $U$ and $V$ are linear spaces.
More generally,
let $X_1,X_2, \cdots, X_n$ $(n \ge 3)$ be distinct vectors in a finite dimensional Euclidean space $E.$
Putting $G:=\frac{\sum_{i=1}^{n} X_i}{n}$, the centroid of the $n$ distinct vectors,
then we get the following identity by a simple direct calculation and the definition of the norm
\begin{eqnarray*}
\sum_{1 \le i<j \le n} \|\overrightarrow{X_i X_j}\|^2
=n \sum_{i=1}^{n} \|\overrightarrow{X_i G}\|^2
\end{eqnarray*}
which is equivalent to the equation
\begin{eqnarray}\label{qu-2}
n \sum_{1 \le i<j \le n} \left\|X_i -X_j \right\|^2
  = \sum_{i=1}^{n}\left\|\sum_{j =1}^n X_j -n X_i\right\|^2
\end{eqnarray}
for any distinct vectors $X_1,X_2, \cdots, X_n$.
Employing the above equality (\ref{qu-2}), we introduce the new functional equation,
\begin{eqnarray}\label{fe3}
 n \sum_{1 \le i<j \le n} Q\left( x_i-x_j \right)=\sum_{i=1}^{n}Q\left(\sum_{j =1}^n x_j -n x_i\right)
\end{eqnarray}
for a mapping $Q:U \rightarrow V$ and for all vectors $x_1,\cdots,x_n \in U.$

%=======================   Quasi-Banach spaces          ===========================================

We recall some basic facts concerning quasi-Banach spaces and some
preliminary results.

\begin{defn} {\rm (\cite{bl00, ro84})} Let $X$
 be a linear space. A {\it quasi-norm } $\|\cdot \|$ is a real-valued function on $X$ satisfying the following:

 (1)  $\|x\|\ge 0$ for all $x \in X$ and $\|x\|=0$ if and
only if $x=0$.

 (2) $\|\lambda x\| = |\lambda| \cdot \|x\|$ for all $\lambda \in \mathbb{R}$ and all $x\in X$.

 (3) There is a constant $K $ such that $\|x+y\| \le K(\|x\| + \|y\|)$ for all $x, y \in X$.

The smallest possible $K$ is called the {\it modulus of concavity } of $\| \cdot \|$.
The pair $(X, \|\cdot \|)$ is called a {\it quasi-normed space} if
$\|\cdot \|$ is a quasi-norm on  $X$.
A {\it quasi-Banach space} is a complete quasi-normed space.
A quasi-norm $\|\cdot \|$ is called a {\it $p$-norm} $(0 < p \le 1)$ if
$$\|x+y\|^p \le \|x\|^p + \|y\|^p $$
for all $x, y \in X$.  In this case, a quasi-Banach space is
called a {\it $p$-Banach space}.
\end{defn}

Clearly, $p$-norms are continuous, and in fact, if $\|\cdot \|$ is a $p$-norm on $X$,
then the formula $d(x, y): = \|x-y\|^p$ defines an
translation invariant metric for $X$ and $\|\cdot \|^p$ is a $p$-homogeneous $F$-norm.
The Aoki--Rolewicz theorem
\cite{bl00,ro84} guarantees that each quasi-norm is equivalent
to some $p$-norm for some $0 < p \le 1$.
Concerning the Ulam stability problem for functional equations,
C. S\'{a}nchez \cite{sa02} and J. Tabor  \cite{ta04} have investigated a version of
the Hyers-Rassias-Gajda theorem (see \cite{ga91, ra78}) for approximate additive mappings in
quasi-Banach spaces.

%=====================================================================

In this paper,
we are going to find the general solution of (\ref{fe3}) for any fixed positive integer $n \ge 3$
in the class of mappings between real vector spaces.
Furthermore,
concerning the stability problem of Ulam for the functional equation (\ref{fe3})
we are going to investigate the generalized Hyers-Ulam-Rassias stability problem
for approximate mappings in quasi-Banach modules and $p$-Banach modules over a $C^*$-algebra.
Thus we generalize the stability results of the quadratic functional equation (\ref{fe3})
in Banach spaces.

% ----------------------------------------------------------------

 \section{Solution of FE. (\ref{fe3})}

First of all, we find out the general solution of (\ref{fe2})
in the class of mappings between real vector spaces.

\begin{lem}\label{aqf0}
Let $U$ and $V$ be real vector spaces.
A mapping $Q:U \rightarrow V$ satisfies the functional equation (\ref{fe2})
if and only if the mapping $Q:U \rightarrow V$ is quadratic.
%Thus in this case there exists a symmetric biadditive mapping $B:U \times U \rightarrow V$
%such that $Q(x)=B(x,x)$ for all $x \in U.$
\end{lem}

{\bf\it Proof.}
It is easy to see that the equation (\ref{fe1}) implies the functional equation (\ref{fe2}).
Now let $Q$ satisfy the equation (\ref{fe2}).
Putting $y,z:=0$ in (\ref{fe2}) yields $Q(2x)=4Q(x)$ for all $x \in U.$
By setting $z:=0$ in (\ref{fe2}), we see
\begin{eqnarray}\label{aqf1-2}
 Q(x-2y)+Q(2x-y)+Q(x+y)
 =3Q(x-y)+3Q(x)+3Q(y)
\end{eqnarray}
for all $x,y \in U.$
In turn, substituting $-y$ for $y$ in (\ref{aqf1-2}) and then adding the resulting equation to (\ref{aqf1-2}),
one obtains
\begin{eqnarray}\label{aqf1-3}
&& Q(2x+y)+Q(2x-y)+Q(x+2y)+Q(x-2y)\\\nonumber
&&\qquad =2Q(x-y)+2Q(x+y)+6Q(x)+6Q(y)
\end{eqnarray}
for any $x,y \in U.$
Letting $z:=-y$ in (\ref{fe2}), we obtain
\begin{eqnarray}\label{aqf1-4}
 Q(x+3y)+Q(x-3y)+4Q(x)
 =3Q(x-y)+3Q(x+y)+12Q(y)
\end{eqnarray}
for all $x,y \in U.$
Replacing $x$ by  $2x$ in (\ref{aqf1-4}), we get
\begin{eqnarray}\label{aqf1-5}
 &&Q(2x+3y)+Q(2x-3y)+16Q(x)\\\nonumber
&&\qquad   =3Q(2x-y)+3Q(2x+y)+12Q(y)
\end{eqnarray}
for any $x,y \in U.$
Now we substitute $z:=2y$ in (\ref{fe2}) to get
\begin{eqnarray}\label{aqf1-6}
 Q(x-3y)+Q(2x-3y)+Q(x)
 =3Q(x-y)+3Q(y)+3Q(x-2y)
\end{eqnarray}
for any $x,y \in U.$
Replacing $y$ by $-y$ in (\ref{aqf1-6}) and then adding (\ref{aqf1-6}) to the resulting expression,
we obtain
\begin{eqnarray*}
&& Q(x+3y)+Q(x-3y)+Q(2x+3y)+Q(2x-3y)+2Q(x)\\\nonumber
&&\qquad =3Q(x+y)+3Q(x-y)+3Q(x+2y)+3Q(x-2y)+6Q(y),
\end{eqnarray*}
which is rearranged in the following way by (\ref{aqf1-5})
\begin{eqnarray}\label{aqf1-7}
&& Q(x+3y)+Q(x-3y)+3Q(2x+y)+3Q(2x-y)+6Q(y)\\\nonumber
&&\qquad =3Q(x+y)+3Q(x-y)+3Q(x+2y)+3Q(x-2y)+14Q(x)
\end{eqnarray}
for any $x,y \in U.$
Now subtracting (\ref{aqf1-4}) from the equation (\ref{aqf1-7}) and then dividing it by $3$,
we have
\begin{eqnarray}\label{aqf1-8}
 Q(2x+y)+Q(2x-y)+6Q(y)=Q(x+2y)+Q(x-2y)+6Q(x)
\end{eqnarray}
for any $x,y \in U.$
Again we add (\ref{aqf1-3}) to (\ref{aqf1-8}) and then divide the resulting expression by $2$
to obtain
\begin{eqnarray}\label{aqf1-9}
 Q(2x+y)+Q(2x-y)=Q(x+y)+Q(x-y)+6Q(x),
\end{eqnarray}
which is equivalent to the original quadratic functional equation
\begin{eqnarray*}
 Q(x+y)+Q(x-y)=2Q(x)+2Q(y)
\end{eqnarray*}
for any $x,y \in U$ \cite[Theorem 2.1]{ck02}.
\hfill$\square$

\bigskip

%==============================  Lemma 2.2 Solution   ========================

\begin{lem}\label{aqf1}
Assume that a mapping $Q:U \rightarrow V$ satisfies the functional equation (\ref{fe3}).
Then $Q$ is even and
\begin{eqnarray}\label{qa1}
Q((n-1)^k x)=(n-1)^{2k} Q(x)
\end{eqnarray} for any vector $x \in U.$
\end{lem}

{\bf\it Proof.}
By setting $x_i:=0$ for all $i=1,\cdots,n$ in the equation (\ref{fe3}), we see $Q(0)=0.$
Putting $x_1 =x$ and $x_i:=0$ for all $i=2,\cdots,n$ in (\ref{fe3}), we get $Q(-(n-1)x)=(n-1)^2Q(x)$ for all $x \in U.$
Substituting $x_i:=x$ for all $i=1,\cdots,n-1$ and $x_n:=0$  in (\ref{fe3}),
one obtains
\begin{eqnarray*}
n(n-1)Q(x)=(n-1)Q(-x)+Q((n-1)x)=n(n-1)Q(-x),
\end{eqnarray*}
which shows that $Q$ is even, and hence $Q((n-1)x)=(n-1)^2 Q(x)$
 for all $x \in U.$
Therefore we get the desired conclusion  by induction on $k.$
\hfill$\square$

\medskip

To find the general solution of (\ref{fe3}), we need to prove the following lemma above all.

%==============================  Lemma 2.3 Solution   ========================

\begin{lem}\label{aqf1-1}
Let $U$ and $V$ be real vector spaces.
For each integer $a$ with $|a| \ne 1,$ a mapping $Q:U \rightarrow V$ satisfies the functional equation
\begin{eqnarray}\label{fe3-0}
&&Q(ax+y)+Q(x+ay)+(a-1)Q(x-y)\\\nonumber
&&\qquad  =(a+1)Q(x+y)+(a^2 -1)[Q(x)+Q(y)]
\end{eqnarray}
 for all $x, y \in U$
if and only if a mapping $Q:U \rightarrow V$ is quadratic.
\end{lem}

{\bf\it Proof.}
Let $Q$ satisfy the equation (\ref{fe3-0}).
It follows easily that $Q$ is even, $Q(ax)=a^2Q(x)$ and $Q(0)=0.$
For $a=0,$ the equation (\ref{fe3-0}) reduces to the equation (\ref{fe1}).
For any negative integer $a< -1$, by considering $a$ as $-a$ and applying the evenness of $Q,$
we need to prove the lemma for the case $a> 1$ without loss of generality.
Now we claim that if $Q$ satisfies the equation (\ref{fe3-0}),
then $Q$ also satisfies (\ref{fe1}) by induction on positive integers $a > 1.$
For $a=2,$ the equation (\ref{fe3-0}) reduces to
\begin{eqnarray}\label{quad-0}
Q(2x+y)+Q(x+2y)+Q(x-y)=3Q(x+y)+3Q(x)+3Q(y),
\end{eqnarray}
which is exactly the equation (\ref{fe2}), and hence it is equivalent to (\ref{fe1}) by Lemma \ref{aqf0}.
Assume that the equation (\ref{fe3-0}) implies the equation (\ref{fe1}) for all $a$ with $a:=2,\cdots,a.$
We are to show that if $Q$ satisfies the equation (\ref{fe3-0}) for $a+1,$ then $Q$ is quadratic in the sequel.
Letting $y:=x+y$ in (\ref{fe3-0}), we obtain
\begin{eqnarray}\label{quad-1}
&&Q((a+1)x+y)+Q((a+1)x+ay)+(a-1)Q(y) \\\nonumber
&& \qquad =(a+1)Q(2x+y)+(a^2 -1) [Q(x)+Q(x+y)]
\end{eqnarray}
 for all $x, y \in U.$
Interchanging $x$ with $y$ in (\ref{quad-1}) and after that adding it to (\ref{quad-1}), we have
\begin{eqnarray}\label{quad-2}
&&Q((a+1)x+y)+Q(x+(a+1)y)+Q((a+1)x+ay)+Q(ax+(a+1)y)\\\nonumber
&&  =(2a^2+3a+1)Q(x+y)+(a^2 +2a+3) [Q(x)+Q(y)]-(a+1)Q(x-y)
\end{eqnarray}
for all $x, y \in U.$
Letting $y:=-x+y$ in (\ref{fe3-0}), we obtain
\begin{eqnarray}\label{quad-3}
&&Q((a-1)x+y)+Q((a-1)x-ay)+(a-1)Q(2x-y)  \\\nonumber
&& \qquad =(a+1)Q(y)+(a^2 -1) [Q(x)+Q(x-y)]
\end{eqnarray}
 for all $x, y \in U.$
Exchanging $x$ and $y$ in (\ref{quad-3}) and after that
adding the resulting equation and (\ref{quad-3}), one has by induction
\begin{eqnarray}\label{quad-4}
&&Q((a-1)x+ay)+Q(ax+(a-1)y)+Q(x-y)  \\\nonumber
&& \qquad =(2a^2-2a-1)Q(x+y)+3[Q(x)+Q(y)]
\end{eqnarray}
 for all $x, y \in U.$
We observe from these inequalities that
$Q( \lambda x)=\lambda ^2 Q(x)$ for $\lambda :=a+1,a-1, 2, 2a-1,$ and for all $x \in U.$
Replacing $y$ by $ay$ in (\ref{quad-4}) and switching $x$ with $y$ in the resulting equation, and then
adding two equations side by side, we obtain by inductive assumption that
\begin{eqnarray}\label{quad-5}
&&Q((a-1)x+a^2 y)+Q(a^2 x+(a-1)y)+(a^3-2a^2+2a+2)Q(x-y)  \\\nonumber
&& \ =(a^3-2a-2)Q(x+y)+(a^4-a^2+2a+5)[Q(x)+Q(y)]
\end{eqnarray}
 for all $x, y \in U.$
Now we substitute $y:=ay-x$ in (\ref{fe3-0}) to get
\begin{eqnarray}\label{quad-6}
&&Q((a-1)x+ay)+Q((a-1)x-a^2 y)+(a-1)Q(2x-ay)  \\\nonumber
&& \qquad =a^2(a+1)Q(y)+(a^2-1)[Q(x)+Q(x-ay)]
\end{eqnarray}
 for all $x, y \in U.$
Switching $x$ with $y$ in (\ref{quad-6}), and then
adding two equations side by side, we obtain by virtue of (\ref{quad-4}), (\ref{quad-5}) and (\ref{fe3-0}) that
\begin{eqnarray}\label{quad-7}
&&(a-1)[Q(2x-ay) +Q(ax-2y)]+(3a^2-5a-2)Q(x+y)  \\\nonumber
&& \qquad =(a^3+a^2-2a-8)[Q(x)+Q(y)]+(a^2+a+2)Q(x-y)
\end{eqnarray}
holds for all $x, y \in U.$
Note from (\ref{quad-7}) that $Q(\lambda x)=\lambda^2 Q(x)$ for $\lambda :=a-2,a+2,$ and for all $x \in U$
Now substituting $x$ for $2x$ in (\ref{fe3-0}) yields
\begin{eqnarray}\label{quad-8}
&&Q(2ax+y) +Q(2x+ay)+(a-1)Q(2x-y)  \\\nonumber
&& \qquad =(a+1)Q(2x+y)+(a^2-1)[4Q(x)+Q(y)]
\end{eqnarray}
for all $x, y \in U.$
Exchanging $x$ and $y$ in (\ref{quad-8}) and then
adding two equations side by side, one obtains by  (\ref{quad-0}), (\ref{quad-7})
\begin{eqnarray}\label{quad-9}
&&(a-1)[Q(2ax+y) +Q(x+2ay)]+(a^2-a+4)Q(x-y)  \\\nonumber
&& \qquad =(4a^3-6a^2+3a+7)[Q(x)+Q(y)]+(3a^2-3a-4)Q(x+y)
\end{eqnarray}
for all $x, y \in U.$ We remark that $Q(\lambda x)=\lambda^2 Q(x)$ for $\lambda :=2a-1,2a+1,$ and for all $x \in U$

Now, let's transform by variables like as $x:=2ax+y$ and $y:=x+2ay$ in (\ref{quad-4}).
Then we can rewrite the equation (\ref{quad-4}) in the form
\begin{eqnarray}\label{quad-10}
&&(2a-1)^2[Q((a+1)x+ay)+Q(ax+(a+1)y)]+(2a-1)^2 Q(x-y)  \\\nonumber
&& \qquad =(2a^2-2a-1)(2a+1)^2Q(x+y)+3[Q(2ax+y)+Q(x+2ay)]
\end{eqnarray}
for all $x, y \in U.$
Multiplying both sides of (\ref{quad-10}) by $(a-1)$ and applying (\ref{quad-9}) to the resulting expression,
we get
\begin{eqnarray}\label{quad-11}
&&(a-1)(2a-1)^2[Q((a+1)x+ay)+Q(ax+(a+1)y)] \\\nonumber
&& \quad +(a-1)(2a-1)^2 Q(x-y)  \\\nonumber
&&  =(8a^5-8a^4-10a^3+13a^2-4a-11)Q(x+y)\\\nonumber
&& \quad  +3(4a^3-6a^2+3a+7)[Q(x)+Q(y)]-(4a^3-5a^2+2a+11)Q(x-y)
\end{eqnarray}
for all $x, y \in U.$
Multiplying $(a-1)(2a-1)^2$ on both sides of (\ref{quad-2}) and
applying (\ref{quad-11}) to the resulting expression, we get finally
\begin{eqnarray*}\label{quad-12}
&&(a-1)(2a-1)^2[Q((a+1)x+y)+Q(x+(a+1)y)] \\\nonumber
&&  =(4a^4-8a^2+6a+10)Q(x+y)-(4a^4-8a^3+2a^2+2a-12)Q(x-y)\\\nonumber
&& \quad  +(4a^5-11a^3+3a^2+4a-24)[Q(x)+Q(y)],
\end{eqnarray*}
which can be written in the form
\begin{eqnarray}\label{quad-13}
&&(a-1)(2a-1)^2\Big[Q((a+1)x+y)+Q(x+(a+1)y) + a Q(x-y) \Big] \\\nonumber
&&  =(a-1)(2a-1)^2\Big[(a+2)Q(x+y)+((a+1)^2-1)[Q(x)+Q(y)]\Big]\\\nonumber
&& \quad +(3a^2-3a+12)[Q(x+y)+Q(x-y)-2Q(x)-2Q(y)]
\end{eqnarray}
for all $x, y \in U.$
Since $Q$ satisfies the equation (\ref{fe3-0}) for $a+1,$
the last equation reduces to $Q(x+y)+Q(x-y)-2Q(x)-2Q(y)=0$.
Consequently, we have proved that if $Q$ satisfies the equation (\ref{fe3-0}) for $a+1,$
then $Q$ satisfies the equation (\ref{fe1}).
Therefore, by induction argument the equation (\ref{fe3-0}) implies the equation (\ref{fe1})
for each positive integer $a>1$.

Conversely, it is obvious that the equation (\ref{fe1}) implies the functional equation (\ref{fe3-0}).
This completes the proof.
\hfill$\square$

\medskip

%==============================  Lemma 2.3 Solution   ========================

\begin{thm}\label{aqf2}
Let $U$ and $V$ be real vector spaces.
A mapping $Q:U \rightarrow V$ satisfies the functional equation (\ref{fe3}) for each positive integer $n>2$
if and only if a mapping $Q:U \rightarrow V$ satisfies the functional equation (\ref{fe1}).
Thus there exists a symmetric biadditive mapping $B:U \times U \rightarrow V$
such that $Q(x)=B(x,x)$ for all $x \in U.$
\end{thm}

{\bf\it Proof.}
It is easy to see that the equation (\ref{fe1}) implies the functional equation (\ref{fe3}).
Conversely, let $Q$ satisfy the equation (\ref{fe3}).
Putting $x_1 :=x$, $x_2 :=y$ and $x_i:=0$ for all $i=3,\cdots,n$ in (\ref{fe3}),
we get
\begin{eqnarray}\label{aqf2-10}
&&Q(x-ay)+Q(ax-y)+(a-1)Q(x+y)\\\nonumber
&&\qquad  =(a+1)Q(x-y)+(a^2 -1)[Q(x)+Q(y)]
\end{eqnarray}
 for all $x, y \in U,$ where $a:=n-1$ is a positive integer with $a \ge 2.$
By the previous Lemma \ref{aqf1-1}, the mapping $Q:U \rightarrow V$ satisfies the functional equation (\ref{fe1}).
\hfill$\square$

\medskip

The following result is interesting and useful characterization formulas for
an inner product space among normed linear spaces.

\begin{cor} Let $U$ be a normed linear space. Then the following statements are equivalent:

(a) $U$ is an inner product space.

(b) The norm in $U$ satisfies the condition:
    \begin{eqnarray*}
  && \|ax+y\|^2 +\|x+ay\|^2 +(a-1)\|x-y\|^2 \\
   && \qquad =(a+1)\|x+y\|^2 +(a^2-1)(\|x\|^2 +\|y\|^2)
    \end{eqnarray*}
     for all $x, y \in U$ and for some fixed integer $a$ with $|a| \ne 1$.

(c)
\begin{eqnarray*}
n \sum_{1 \le i<j \le n} \|x_i -x_j\|^2 = \sum_{i=1}^{n} \left\|\sum_{j=1}^n x_j -n x_i \right\|^2
 \end{eqnarray*}
for all $x_i (i=1,\cdots, n) \in U$ and a fixed $n>2.$
\end{cor}
{\bf\it Proof.}
The proof is obvious by Lemma \ref{aqf1-1} and Theorem \ref{aqf2}.
The inner product is defined as usual by
\begin{eqnarray*}
 (x,y)&=&1/4\left( \|x+y\|^2 -\|x-y\|^2 +i\|x+iy\|^2-i \|x-iy\|^2 \right), \ \text{and} \\
 (x,y)&=&1/4\left( \|x+y\|^2 -\|x-y\|^2  \right)
\end{eqnarray*}
for the complex and real spaces, respectively.

\hfill$\square$

\medskip

%================================   Section 3    =========================================

 \section{Stability of FE. (\ref{fe3}) in quasi-Banach modules}

Now let $\mathcal{A}$ be a complex $*$-algebra with unit and let $M$ be a left $\mathcal{A}$-module.
Let us call a mapping $Q:M \rightarrow \mathcal{A}$ an {\it $\mathcal{A}$-quadratic mapping} if both relations
$Q(ax)=aQ(x)a^*$ and $Q(x+y)+Q(x-y)=2Q(x)+2Q(y)$ are fulfilled \cite{vu87}.
A mapping $Q:M \rightarrow \mathcal{A}$ is called a {\it generalized $\mathcal{A}$-quadratic mapping} if $Q(ax)=aQ(x)a^*$
for all $x \in M$, and the following identity holds:
\begin{eqnarray*}
Q\left(\sum_{i=1}^n a_ix_i\right)+ \sum_{1\le i<j\le n} a_ia_j Q(x_i -x_j)
=\left(\sum_{i=1}^n a_i\right) \left[\sum_{i=1}^na_i Q(x_i)\right]
\end{eqnarray*}
for all $x_i \in M,$ some fixed $a_i$ in $\Bbb R$  $(i=1, \cdots, n)$
 and at least two of them are nonzero such that $\sum_{i=1}^n a_i \ne 0,$ and a fixed $n \ge 2$ \cite{lin92}.
It was shown that the notion of $\mathcal{A}$-quadratic mapping is equivalent to the notion of
generalized $\mathcal{A}$-quadratic mapping if all spaces are over the complex number field
and a mapping $B:M \times M \rightarrow \mathcal{A}$ is defined in terms of the mapping
$Q$ as
\begin{eqnarray}\label{aqf3-22-0}
 B(x,y)=\frac{1}{4}[Q(x+y)-Q(x-y)+iQ(x+iy)-iQ(x-iy)]
\end{eqnarray}
for all $x,y$ in $M$ \cite{lin92}.
It was indicated in \cite{vu87} that if the relation (\ref{aqf3-22-0}) holds and
$Q$ is an $\mathcal{A}$-quadratic form, then $B$ is an $A$-sesquilinear form and  $Q(x)=B(x,x),$ and vice versa.
Now it follows easily from Theorem \ref{aqf2} that a mapping $Q$ is a generalized $\mathcal{A}$-quadratic mapping if and only if
\begin{eqnarray*}
&&Q(ax)=aQ(x)a^*,\\
&& n \sum_{1\le i<j\le n} Q(x_i-x_j)
=\sum_{i=1}^{n}Q\left(\sum_{j=1}^n x_j -n x_i\right)
\end{eqnarray*}
for all $x$ and $(x_1, \cdots, x_n),$ where $n \ge 3$.

Now we are ready to investigate the generalized Hyers-Ulam-Rassias stability problem
for approximate $\mathcal{A}$-quadratic mappings acting on $\mathcal{U}(\mathcal{A})$
of the equation (\ref{fe3}) in quasi-Banach modules over a $C^*$-algebra.
Let $M_1$ and $M_2$ be quasi-Banach $\mathcal{A}$-bimodules and
let $K \ge 1$ be the  modulus of concavity of $\| \cdot \|$
throughout this section unless we give any specific reference.
Given a mapping $f : M_1 \rightarrow M_2$,
we define a difference $D_uf:M_1^n \rightarrow M_2$ of the equation (\ref{fe3})  as
\begin{eqnarray*}
&&D_uf(x_1,\cdots,x_n)\\
&&\quad :=n\sum_{1\le i<j\le n} f(ux_i-ux_j)
  -\sum_{i=1}^{n}uf\left(\sum_{j =1}^n x_j -n x_i\right)u^*,
\end{eqnarray*}
for all $x_i \in M_1$ and $u \in \mathcal{U}(\mathcal{A}),$ which is called the approximate remainder
of the functional equation (\ref{fe3}) and acts as a perturbation of the equation.

%%%%%%%%%%%%%%%%%%%%%%%%%%%%%%%%%%%55                       %%%%%%%%%%%%%%%%%%%%%%

\begin{thm}\label{aqf5}
Assume that there exists a mapping $\varphi : M_1^n \rightarrow [0, \infty) := \mathbb {R}_+$ for which
 a mapping $f: M_1 \rightarrow M_2$  satisfies the functional inequality
\begin{equation}\label{aqf3-1}
    \|D_uf(x_1,\cdots,x_n)\| \leq \varphi(x_1,\cdots,x_n)
\end{equation}
for all $(x_1,\cdots,x_n) \in M_1^n$ and for all $u \in \mathcal{U}(\mathcal{A}),$
and the following series
\begin{eqnarray}\label{aqf3-1-1}
 \sum _{i=0}^\infty
\frac{K^i \varphi((n-1)^{i}x_1,\cdots,(n-1)^{i}x_n)} {(n-1)^{2i}}<\infty
\end{eqnarray}
for all $(x_1,\cdots,x_n) \in M_1^n$.
If either $f$ is measurable or $f(tx)$ is continuous in $t \in \mathbb{R}$ for each
fixed $x \in M_1,$
then there exists a unique generalized $\mathcal{A}$-quadratic mapping $Q: M_1 \rightarrow M_2$ which satisfies
the equation (\ref{fe3}) and the inequality
\begin{eqnarray}\label{aqf3-1-2}
\left\|f(x)+\frac{(n-1)f(0)}{2}-Q(x)\right\| \le \frac{K}{(n-1)^2}
\sum _{i=0}^\infty
\frac{K^i \Phi((n-1)^{i}x)} {(n-1)^{2i}}
\end{eqnarray}
for  all  $x \in M_1,$ where
\begin{eqnarray*}
 \Phi(x)&:=&\min_{1\le i \le n} \left\{\varphi_{i}(-x) + \frac{|(n^2+1)-(i+1)n|}{n}\tilde{\varphi}(x)\right\},\\
\varphi_{i}(x) &:=& \varphi(0,\cdots,0,\underbrace{x}_{i-th},0,\cdots,0),\ (i=1,\cdots,n),\\
\text{and} \ \ \tilde{\varphi}(x)&:=&\min_{1\le i \le n-1} \left\{\varphi_{i}(x) + \varphi_{i+1}(x)\right\}
\end{eqnarray*}
for  all  $x \in M_1.$
The mapping $Q$ is defined by
\begin{eqnarray*}
 Q(x)=\lim\limits_{k\rightarrow \infty} { f((n-1)^k x)  \over (n-1)^{2k}}
\end{eqnarray*}
for  all  $x \in M_1.$
\end{thm}

%%%%%%%%%%%%%%%%%%%%%%%%%%%   Thm 3.1   %%%%%%%%%%%%%%%%%%%%%%%%%%%%%%%%%%%%%%%%%%%%%%%%%%
%%%%%%%%%%%%%%%%%%%%%%%%%%%  Stability of (\ref{fe3})   %%%%%%%%%%%%%%%%%%%%%%%%%%%%%%%%%%%%%%%%%%%%%%%%

{\bf\it Proof.}
Put $u:=1 \in \mathcal{U}(\mathcal{A})$ in (\ref{aqf3-1}).
Then for each $i=1,\cdots,n-1,$ interchanging $x_i$ for $x$ and $x_j$ for $0$ for all $j \ne i$ in (\ref{aqf3-1})
and then comparing the sequent inequalities, we get
\begin{eqnarray*}
n\|f(x)-f(-x)\| \le \varphi_{i}(x)+\varphi_{i+1}(x)
\end{eqnarray*}
for all $x \in M_1$ and for all $i=1,\cdots,n-1$.
Thus one obtains the approximate even condition of $f$
\begin{eqnarray}\label{aqf3-2}
\|f(x)-f(-x)\| \le \frac{1}{n}\tilde{\varphi}(x),
\quad \tilde{\varphi}(x):=\min_{1\le i \le n-1} \left\{\varphi_{i}(x) + \varphi_{i+1}(x)\right\}
\end{eqnarray}
for all $x \in M_1.$
For each $i=1,\cdots,n,$ replacing $x_i$ by $-x$ and $x_j$ by $0$ for all $j \ne i$ we observe that
\begin{eqnarray*}
&&\left\|(i-1)n f(x) +[n^2 -(i+1)n +1]f(-x)
+n {n-1 \choose  2} f(0)-f((n-1)x)\right\| \\\nonumber
&&\quad \le \varphi_{i}(-x)
\end{eqnarray*}
for all $x \in M_1.$
Associating the last inequality with (\ref{aqf3-2}), we obtain
\begin{eqnarray*}
&&\left\|(n-1)^2f(x)
+n {n-1 \choose  2}
  f(0)-f((n-1)x)\right\| \\\nonumber
&&\qquad \le \varphi_{i}(-x) + \frac{|(n^2+1)-(i+1)n|}{n}\tilde{\varphi}(x)
\end{eqnarray*}
for all $x \in M_1$ and for all $i=1,\cdots,n$.
Hence one has the following inequality
\begin{eqnarray}\label{aqf3-2-1}
&&\left\|(n-1)^2f(x)
+n {n-1 \choose  2}
  f(0)-f((n-1)x)\right\|  \le \Phi(x)
\end{eqnarray}
for all $x \in M_1.$
Define a mapping $g : M_1 \rightarrow M_2$ by
 $g(x):= f(x)+\frac{(n-1)f(0)}{2}$ for all $x \in M_1$.
Then it follows from (\ref{aqf3-2-1}) that
\begin{eqnarray}\label{aqf3-3}
&&\left\|g(x)-\frac{g((n-1)x)}{(n-1)^2}\right\|  \le \frac{1}{(n-1)^2}\Phi(x)
\end{eqnarray}
for all $x \in M_1,$
from which
we obtain by applying a standard procedure of the induction argument on $m$ that
\begin{eqnarray}\label{aqf3-4-0}
\left\|g( x)  -{g((n-1)^{m} x)\over (n-1)^{2m}} \right\|
&\le& \frac{K}{(n-1)^{2}} \sum_{i=0}^{m-2} \left(\frac{K}{(n-1)^{2}}\right)^i \Phi((n-1)^ix) \\\nonumber
&+& \frac{1}{(n-1)^{2}}\left(\frac{K}{(n-1)^{2}}\right)^{m-1} \Phi((n-1)^{m-1}x )
\end{eqnarray}
for all $x \in M_1$ and all $m \ge 1,$ which is considered to be (\ref{aqf3-3}) for $m=1.$
In fact, we figure out by the inequality (\ref{aqf3-3}),
\begin{eqnarray*}
&&\left\|g( x)  -{g((n-1)^{m+1} x)\over (n-1)^{2(m+1)}} \right\|\\\nonumber
&&\le K\left\|g( x)  -{g((n-1) x)\over (n-1)^{2}} \right\|
    + K \left\|{g((n-1) x)\over (n-1)^{2}} -{g((n-1)^{m+1} x)\over (n-1)^{2(m+1)}} \right\| \\\nonumber
&&\le \frac{K}{(n-1)^2}\Phi(x)
    + \frac{K}{(n-1)^2}\left\|g((n-1) x) -{g((n-1)^{m+1}x) \over (n-1)^{2m}} \right\|,
\end{eqnarray*}
which, in accordance with inductive assumption, yields (\ref{aqf3-4-0}) for $m+1.$
Thus one obtains that for all nonnegative integers $m, l$ with $m>l$
\begin{eqnarray}\label{aqf3-4}
&&\left\|{g((n-1)^l x) \over (n-1)^{2l}} -{g((n-1)^{m} x)\over (n-1)^{2m}} \right\| \\\nonumber
&&= {1 \over (n-1)^{2l}} \left\|g( (n-1)^l x)  -{g((n-1)^{m-l}\cdot (n-1)^{l}x)\over (n-1)^{2(m-l)}} \right\| \\\nonumber
&&\le \frac{K}{(n-1)^{2l+2}} \sum_{i=0}^{m-l-2} \frac{K^i \Phi((n-1)^{l+i}x)}{ (n-1)^{2i}}
+\frac{1}{(n-1)^{2l+2}}  \frac{K^{m-l-1} \Phi((n-1)^{m-1}x) }{ (n-1)^{2(m-l-1)}} \\\nonumber
&&\le \frac{K}{K^l(n-1)^{2}} \sum_{i=l}^{m-2} \frac{K^i \Phi((n-1)^{i}x)}{ (n-1)^{2i}}
+\frac{1}{K^l(n-1)^{2}}  \frac{K^{m-1} \Phi((n-1)^{m-1}x) }{ (n-1)^{2(m-1)}},
\end{eqnarray}
which tends to zero by (\ref{aqf3-1-1}) as $l \rightarrow \infty.$
Hence the sequence $\Big\{ {g((n-1)^m x ) \over (n-1)^{2m}} \Big\}_{m\in \Bbb N}$ is a Cauchy sequence
for any $x\in M_1,$ and so it converges by the completeness of $M_2$.
Therefore we can define a mapping $Q: M_1 \rightarrow M_2$ by
\begin{eqnarray*}
 Q(x)=\lim\limits_{m\rightarrow \infty} { g((n-1)^m x) \over (n-1)^{2m}}
 =\lim\limits_{m\rightarrow \infty} { f((n-1)^m x)  \over (n-1)^{2m}}
\end{eqnarray*}
for all $x \in M_1$.
Taking the limit as $m \rightarrow \infty$ in (\ref{aqf3-4-0}), we
obtain the desired inequality (\ref{aqf3-1-2}).
Exchanging $(x_1,\cdots,x_n)$ for $((n-1)^mx_1,\cdots,(n-1)^mx_n)$ in (\ref{aqf3-1}) and dividing
both sides by $(n-1)^{2m}$, we have
\begin{eqnarray}\label{aqf3-8}
&&\|D_1Q(x_1,\cdots,x_n)\| \\\nonumber
&&\quad =\lim\limits_{m\rightarrow \infty} \frac{1}{(n-1)^{2m}}\|Df((n-1)^mx_1,\cdots,(n-1)^mx_n)\|\\\nonumber
&&\quad \le \lim\limits_{m\rightarrow \infty} \frac{K^m}{(n-1)^{2m}}\varphi((n -1)^mx_1,\cdots,(n-1)^mx_n)\\\nonumber
&&\quad =0.
\end{eqnarray}
Therefore the mapping $Q$ satisfies the equation (\ref{fe3}) and hence $Q$ is quadratic.

To prove the uniqueness, let $Q^{\prime}$ be another quadratic
mapping satisfying (\ref{aqf3-1-2}).
Then we get by Lemma \ref{aqf1} that
$Q^{\prime}((n-1)^mx)=(n-1)^{2m}Q^{\prime}(x)$ for all $x\in M_1$ and all $m \in \mathbb{N}$.
Thus we have
\begin{eqnarray*}
&&\|Q(x)-Q^{\prime}(x)\| \\
&&  \le {1 \over (n-1)^{2m}} \Big\{ K\Big\|Q((n-1)^mx) -f((n-1)^m x)-\frac{(n-1)f(0)}{2} \Big\| \\
&&\qquad\qquad \qquad\qquad  + K\Big\|f((n-1)^m x)+\frac{(n-1)f(0)}{2}-Q^{\prime}((n-1)^m x) \Big\| \Big\}\\\
&& \leq {2K^2 \over K^m(n-1)^{2}}\sum_{i=0}^{\infty}
\frac{ K^{m+i} \Phi((n-1)^{k+i}x )} {(n-1)^{2(m+i)} }
 \end{eqnarray*}
for all $x\in M_1$. Taking the limit as $m \rightarrow \infty$, then we
conclude that $Q(x)=Q^{\prime}(x)$ for all $x\in M_1$.

Finally, we show that the quadratic mapping $Q$ is $\mathcal{A}$-quadratic.
Under the assumption that either $f$ is measurable or
$f(tx)$ is continuous in $t \in \mathbb{R}$ for each
fixed $x \in M_1$,  the quadratic mapping $Q$ satisfies
$Q(tx) = t^2Q(x)$  for all $ x\in M_1$ and all $t \in \mathbb{R}$ by the same reasoning as the proof of \cite{cz92}.
That is, $Q$ is $\mathbb{ R}$-quadratic.
Putting $x_1:=-(n-1)^m x$ and $x_i:=0$ for all $i=2,\cdots,n$ in (\ref{aqf3-1})
and dividing the resulting inequality by $(n-1)^{2m},$
\begin{eqnarray*}
&&\frac{1}{(n-1)^{2m}}\Big\|n(n-1)f(-(n-1)^m ux)+n  {n-1 \choose  2}  f(0) \\
&&\qquad \qquad \qquad -uf((n-1)^{m+1}x)u^* -(n-1)u f(-(n-1)^m x)u^* \Big\| \\
&&\qquad\qquad \le \frac{K^m}{(n-1)^{2m}}\varphi(-(n-1)^mx,0,\cdots,0).
\end{eqnarray*}
Taking the limit as $m \rightarrow \infty$ and using the evenness of $Q$,
we see that $Q(u x)=uQ(x)u^*$ for all $x \in M_1$ and for each $u \in \mathcal{U}(\mathcal{A}).$
The last relation is also true for $u=0.$
Now let $a$ be a nonzero element in $ \mathcal{A}$ and $L$ a positive integer greater than $4|a|.$
Then we have $|\frac{a}{L}|<\frac{1}{4}<1-\frac{2}{3}.$
By \cite[Theorem 1]{kp85}, there exist three elements $u_1,u_2,u_3 \in \mathcal{U}(\mathcal{A})$ such that
$3\frac{a}{L}=u_1 +u_2 +u_3.$
Thus we calculate in conjunction with \cite[Lemma 2.1]{jk05} that
\begin{eqnarray*}
 Q(ax)  &=& Q \left( \frac{L}{3} 3\frac{a}{L}x \right)= \left(\frac{L}{3}\right)^2 Q(u_1x+u_2x+u_3x)  \\
&=& \left(\frac{L}{3}\right)^2 B(u_1x+u_2x+u_3x,u_1x+u_2x+u_3x) \\
&=& \left(\frac{L}{3}\right)^2 (u_1+u_2+u_3)B(x,x)(u_1^*+u_2^*+u_3^*) \\
&=& \left(\frac{L}{3}\right)^2 3\frac{a}{L} Q(x) 3\frac{a^*}{L}= aQ(x)a^*
\end{eqnarray*}
for all $a \in \mathcal{A} (a \ne 0)$ and for all $x \in M_1.$
 So the unique $\mathbb{R}$-quadratic mapping
$Q$ is also generalized $\mathcal{A}$-quadratic, as desired. This completes the proof.
\hfill$\square$

%%%%%%%%%%%%%%%%%%%%%%%%%%%%%%%%%%%   Theorem 3.2                   %%%%%%%%%%%%%%%%%%%%%%

\begin{thm}\label{aqf4}
 Assume that the approximate remainder $D_uf$ of a mapping
 $f : M_1 \rightarrow M_2$ satisfies the functional inequality (\ref{aqf3-1})
 for all $(x_1,\cdots,x_n) \in M_1^n$ and all $u \in \mathcal{U}(\mathcal{A}),$
 and that the following series
\begin{eqnarray}\label{aqf3-23}
\sum _{i=1}^{\infty} K^i(n-1)^{2i} \varphi\left({x_1 \over (n-1)^{i}},\cdots,{x_n \over (n-1)^{i}}\right)
\end{eqnarray}
converges for all $(x_1, \cdots, x_n) \in M_1^n$.
If either $f$ is measurable or $f(tx)$ is continuous in $t \in \mathbb{R}$ for each
fixed $x \in M_1,$ then there exists a unique generalized $\mathcal{A}$-quadratic mapping
$Q: M_1 \rightarrow M_2$ which satisfies the equation (\ref{fe3})
and the inequality
\begin{eqnarray}\label{aqf3-24}
\|f(x)-Q(x)\|
\leq {1 \over (n-1)^2} \sum_{i=1}^{\infty} K^i(n-1)^{2i} \Phi \left(\frac{x}{(n-1)^i} \right)
\end{eqnarray}
 for all $x \in M_1,$ where $\Phi$ is defined as in Theorem \ref{aqf5}.
 The mapping $Q$ is defined by
\begin{eqnarray*}
 Q(x)=\lim\limits_{m \rightarrow \infty}  (n-1)^{2m} f\Big( \frac{x}{(n-1)^m}  \Big)
\end{eqnarray*}
 for all $x \in M_1.$
\end{thm}

%%%%%%%%%%%%%%%%%%%%%%%%%%%   Thm 3.2   %%%%%%%%%%%%%%%%%%%%%%%%%%%%%%%%%%%%%%%%%%%%%%%%%%
%%%%%%%%%%%%%%%%%%%%%%%%%%%             %%%%%%%%%%%%%%%%%%%%%%%%%%%%%%%%%%%%%%%%%%%%%%%%

{\bf\it Proof.} We use the same notations as those of Theorem \ref{aqf5}.
We observe that $\varphi(0,\cdots,0)=0$ by the convergence (\ref{aqf3-23}), and thus we have
$f(0)=0$ by setting $x_i:=0$ in (\ref{aqf3-1}) for all $i=1,\cdots,n$.
Now we get by (\ref{aqf3-3})
\begin{eqnarray*}\label{aqf3-11}
\left\|f(x) - (n-1)^2f\left(\frac{x}{n-1} \right)\right\| \le \Phi\left(\frac{x}{n-1}\right),\quad x \in M_1,
\end{eqnarray*}
which yields by induction
\begin{eqnarray}\label{aqf3-12}
&&\left\|f(x)-(n-1)^{2m} f\left( \frac{x}{(n-1)^m} \right) \right\| \\\nonumber
&&\le  \frac{1}{(n-1)^2}\sum_{i=1}^{m-1} K^i(n-1)^{2i} \Phi\left(\frac{x}{(n-1)^i}\right)
  + \frac{K^m(n-1)^{2m}}{K (n-1)^2}  \Phi\left(\frac{x}{(n-1)^m}\right)
\end{eqnarray}
for all $x \in M_1$ and all integers $m >1$.

It follows by (\ref{aqf3-23}) that the sequence
$\left\{ (n-1)^{2m} f\left({x \over (n-1)^m}\right) \right\}_{m\in \Bbb N}$ is a Cauchy
sequence for any $x\in M_1.$ Since $M_2$ is complete, we may
define a mapping $Q:M_1 \rightarrow M_2$ by
\begin{eqnarray*}
Q(x)=\lim\limits_{m\rightarrow \infty} (n-1)^{2m} f\left({x \over (n-1)^m}\right),\ x \in M_1.
\end{eqnarray*}

The rest of the proof goes through by
 the same way as that of Theorem \ref{aqf5}.
 This completes the proof.
 \hfill$\square$

\medskip

From the main Theorem \ref{aqf5} and Theorem \ref{aqf4} we obtain the following
corollary concerning the stability of the equation (\ref{fe3}).

%%%%%%%%%%%%%%%%%%%%%%%%%%%%%%  Corollary    %%%%%%%%%%%%%%%%%%%%%%%%%%%%
%%%%%%%%%%%%%%%%%%%%%%%%%%%%%%               %%%%%%%%%%%%%%%%%%%%%%%%%%%

\begin{cor}\label{aqfcor1}
Let $r,\varepsilon$ be positive real numbers with  $r-2 < -\log_{n-1}K$ or $r-2 > \log_{n-1}K.$
Assume that a mapping $f : M_1 \rightarrow M_2$ satisfies the
inequality
\begin{equation}\label{aqf3-21}
    \|D_uf(x_1,\cdots,x_n)\| \leq  \varepsilon \sum_{i=1}^{n} \|x_i\|^r
\end{equation}
for all $x_i \in M_1$ and for all $u \in \mathcal{U}(\mathcal{A}).$
Then there exists a unique generalized $\mathcal{A}$-quadratic mapping $Q:M_1 \rightarrow M_2$
which satisfies the equation (\ref{fe3}) and the inequality
\begin{eqnarray*}
\left\|f(x) -Q(x)\right\|
\leq
\left\{
\begin{array}{cc}
 \frac { (n+2)K \varepsilon  \|x\|^r }{n[(n-1)^2-K (n-1)^r] },&  \text{if} \quad  r-2 < -\log_{n-1} K  \\
 \frac { (n+2)K \varepsilon  \|x\|^r }{n[(n-1)^r-K (n-1)^2] },&  \text{if} \quad  r-2 > \log_{n-1} K  \\
\end{array}\right\}
\end{eqnarray*}
for all $x \in M_1$.
\end{cor}

{\bf\it Proof.}
Define $\varphi(x_1,\cdots,x_n):= \varepsilon(\|x_1\|^r + \cdots+ \|x_n\|^r)$
for all $(x_1,\cdots,x_n) \in M_1^n$.
Then we have the conditions
$\tilde{\varphi}(x):=2\varepsilon\|x\|^p$ and $\Phi(x):=\varepsilon\|x\|^p+\frac{2}{n}\varepsilon\|x\|^p$
for all $x \in M_1.$
Applying Theorem \ref{aqf5} and Theorem \ref{aqf4}, we obtain the desired results
according to the cases of $r$.
Exchanging $x_i$ for $0$ in (\ref{aqf3-21}) for all $i=1,\cdots,n$ yields $f(0)=0.$
\hfill$\square$

 \medskip

\begin{prob}
It is an open problem to investigate the stability problem of Ulam
for the case of $K$ and $r$ with $ -\log_{n-1} K \le r-2 \le  \log_{n-1} K$ in Corollary \ref{aqfcor1}.
\end{prob}

%%%%%%%%%%%%%%%%%%%%%%%%%%%%%%  Corollary    %%%%%%%%%%%%%%%%%%%%%%%%%%%%
%%%%%%%%%%%%%%%%%%%%%%%%%%%%%%               %%%%%%%%%%%%%%%%%%%%%%%%%%%

\begin{cor}\label{aqfcor2}
Assume that there exists a nonnegative number $\theta$ for which a mapping $f : M_1 \rightarrow M_2$ satisfies the
inequality
\begin{equation*}
    \|D_uf(x_1,\cdots,x_n)\| \leq  \theta
\end{equation*}
for all $x_i \in M_1$ and for all $u \in \mathcal{U}(\mathcal{A})$.
Then there exists a unique generalized $\mathcal{A}$-quadratic mapping
$Q:M_1 \rightarrow M_2$ which satisfies the equation
(\ref{fe3}) and the inequality
\begin{eqnarray*}
\left\|f(x)+\frac{(n-1)f(0)}{2}-Q(x)\right\|
\leq \frac{(n+2)K\theta}{n[(n-1)^2 -K]} , \ \  \text{if} \ \  K < (n-1)^2
\end{eqnarray*}
for all $x \in M_1$.
\end{cor}

\begin{prob}If $K$ is so large a constant that $(n-1)^2 \le K,$ then we can't guarantee
that the functional equation (\ref{fe3}) is stable on concerning the Ulam stability problem.
So, it is interesting to investigate the stability problem of Ulam
for the case of $n,K$ with $(n-1)^2 \le K$ in Corollary \ref{aqfcor2}.
\end{prob}

%================================   Section 4    =========================================

 \section{Stability of FE. (\ref{fe3}) in $p$-Banach modules}

We now prove the Hyers--Ulam--Rassias stability of the functional equation (\ref{fe3})
in $p$-Banach $\mathcal{A}$-bimodules.

\begin{thm}\label{aqf6}
Let $M_1$ and $M_2$ be $p$-Banach $\mathcal{A}$-bimodules.
Assume that the approximate remainder $D_uf$ of a mapping
 $f : M_1 \rightarrow M_2$ satisfies the functional inequality (\ref{aqf3-1})
 for all $(x_1,\cdots,x_n) \in M_1^n$ and all $u \in \mathcal{U}(\mathcal{A}),$
 and that the following series
\begin{eqnarray*}\label{aqf6-1-1}
 \sum _{i=0}^\infty
\frac{{\varphi((n-1)^{i}x_1,\cdots,(n-1)^{i}x_n)}^p} {(n-1)^{2ip}}<\infty
\end{eqnarray*}
for all $(x_1,\cdots,x_n) \in M_1^n$.
If either $f$ is measurable or $f(tx)$ is continuous in $t \in \mathbb{R}$ for each
fixed $x \in M_1,$
then there exists a unique generalized $\mathcal{A}$-quadratic mapping $Q: M_1 \rightarrow M_2$ which satisfies
the equation (\ref{fe3}) and the inequality
\begin{eqnarray*}\label{aqf6-2}
\left\|f(x)+\frac{(n-1)f(0)}{2}-Q(x)\right\| \le \frac{1}{(n-1)^2}
\left[\sum _{i=0}^\infty
\frac{ {\Phi((n-1)^{i}x)}^p} {(n-1)^{2ip}}\right]^{1/p}
\end{eqnarray*}
for all $x \in M_1,$ where $Q$ and $\Phi$ are defined as in Theorem \ref{aqf5}.
\end{thm}

%%%%%%%%%%%%%%%%%%%%%%%%%%%   Thm 3.8   %%%%%%%%%%%%%%%%%%%%%%%%%%%%%%%%%%%%%%%%%%%%%%%%%%
{\it Proof.}
It follows by the inequality (\ref{aqf3-3}) and the definition of $p$-norm that
\begin{eqnarray*}\label{aqf6-3}
\left\| \frac{g( (n-1)^i x)}{(n-1)^{2i}}  - \frac{g( (n-1)^{i+1} x)}{(n-1)^{2(i+1)}}\right\|^p
  \le   \frac{1}{(n-1)^{2p}} \frac{1}{(n-1)^{2pi}}  \Phi((n-1)^i x)^p ,
\end{eqnarray*}
and so
\begin{eqnarray*}
\left\| \frac{g( (n-1)^l x)}{(n-1)^{2l}} - \frac{g( (n-1)^m x)}{(n-1)^{2m} }\right\|^p
&\le&  \sum_{i=l}^{m-1} \left\| \frac{g( (n-1)^i x)}{(n-1)^{2i}}
- \frac{g( (n-1)^{i+1} x)}{(n-1)^{2(i+1)}}\right\|^p  \\\nonumber
 & \le&   \frac{1}{(n-1)^{2p}} \sum_{i=l}^{m-1}\frac{1}{(n-1)^{2pi}}  \Phi((n-1)^i x)^p
\end{eqnarray*}
for all $x \in M_1$ and all integers $l,m$ with $m > l\ge 0$.
Note that the series $\sum _{i=0}^\infty \frac{ {\Phi((n-1)^{i}x)}^p} {(n-1)^{2ip}}$
converges for all $x \in M_1.$ Thus we obtain the desired results using the similar argument to
Theorem \ref{aqf5}.
\hfill$\square$
%%%%%%%%%%%%%%%%%%%%%%%%%%%%%%%%%%%   Theorem 3.2                   %%%%%%%%%%%%%%%%%%%%%%

\begin{thm}\label{aqf7}
Let $M_1$ and $M_2$ be $p$-Banach $\mathcal{A}$-bimodules.
 Assume that the approximate remainder $D_uf$ of a mapping
 $f : M_1 \rightarrow M_2$ satisfies the functional inequality (\ref{aqf3-1})
 for all $(x_1,\cdots,x_n) \in M_1^n$ and all $u \in \mathcal{U}(\mathcal{A}),$
 and that the following series
\begin{eqnarray*}\label{aqf7-1}
\sum_{i=1}^{\infty} (n-1)^{2ip}
{\varphi\left({x_1 \over (n-1)^{i}},\cdots,{x_n \over (n-1)^{i}}\right)}^{p} < \infty
\end{eqnarray*}
for all $(x_1, \cdots, x_n) \in M_1^n$.
If either $f$ is measurable or $f(tx)$ is continuous in $t \in \mathbb{R}$ for each
fixed $x \in M_1,$ then there exists a unique generalized $\mathcal{A}$-quadratic mapping
$Q: M_1 \rightarrow M_2$ which satisfies the equation (\ref{fe3})
and the inequality
\begin{eqnarray*}
\|f(x)-Q(x)\| \leq {1 \over (n-1)^2}
\left[ \sum_{i=1}^{\infty} (n-1)^{2ip} \Phi \left( \frac{x}{(n-1)^i} \right)^p \right]^{1/p}
\end{eqnarray*}
 for all $x \in M_1,$ where $Q$ and $\Phi$ are defined as in Theorem \ref{aqf4}.
\end{thm}

%%%%%%%%%%%%%%%%%%%%%%%%%%%   Thm 3.9   %%%%%%%%%%%%%%%%%%%%%%%%%%%%%%%%%%%%%%%%%%%%%%%%%%

\begin{cor}\label{aqfcor3}
Let $M_1$ and $M_2$ be $p$-Banach $\mathcal{A}$-bimodules.
Let $r,\varepsilon$ be positive real numbers with  $r \ne 2.$
Assume that a mapping $f : M_1 \rightarrow M_2$ satisfies the
inequality
\begin{equation*}
    \|D_uf(x_1,\cdots,x_n)\| \leq  \varepsilon \sum_{i=1}^{n} \|x_i\|^r
\end{equation*}
for all $x_i \in M_1$ and for all $u \in \mathcal{U}(\mathcal{A}).$
Then there exists a unique generalized $\mathcal{A}$-quadratic mapping $Q:M_1 \rightarrow M_2$
which satisfies the equation (\ref{fe3}) and the inequality
\begin{eqnarray*}
\left\|f(x) -Q(x)\right\|
\leq
\left\{
\begin{array}{cc}
 \frac { (n+2) \varepsilon  \|x\|^r }{n[(n-1)^{2p} -  (n-1)^{rp} ]^{1/p} },&  \text{if} \quad  r < 2 \\
 \frac { (n+2) \varepsilon  \|x\|^r }{n[(n-1)^{rp}-  (n-1)^{2p} ]^{1/p} },&  \text{if} \quad  r > 2  \\
\end{array}\right\}
\end{eqnarray*}
for all $x \in M_1$.
\end{cor}

%=====================================================================================================

\begin{rem}
The result for the case $K=1$ in Theorem \ref{aqf5} (Theorem \ref{aqf4}, Corollary \ref{aqfcor1}, respectively)
is the same as the
result for the case $p=1$ in Theorem \ref{aqf6} (Theorem \ref{aqf7}, Corollary \ref{aqfcor3}, respectively).
\end{rem}

\medskip

Let $M_1$ and $M_2$ be Banach left $A$-modules
and let $\hat{a} :=aa^*,a^*a,$ or $\frac{aa^* +a^*a}{2}$ for each $a \in \mathcal{A}$.
A mapping $Q: M_1 \rightarrow M_2$ is called $\mathcal{A}_{sa}$-quadratic if
$Q(x+y)+Q(x-y)=2Q(x)+2Q(y)$ and
$Q(ax)=\hat{a}Q(x)$ for all $a \in \mathcal{A}$ and all $x,y \in M_1$ \cite{pa02}.
Since two Banach spaces $E_1$ and $E_2$ are considered as Banach modules over $\mathcal{A}:=\Bbb C,$
 the $\mathcal{A}_{sa}$-quadratic mapping $Q: E_1 \rightarrow E_2$ implies $Q(ax)=|a|^2Q(x)$ for all $a \in \Bbb C.$

%%%%%%%%%%%%%%%%%%%%%%%%%%%%%%%%%%%%%         Thm       %%%%%%%%%%%%%%%%%%%%%%%%

\begin{thm}\label{aqf8}
Let $M_1$ and $M_2$ be quasi-Banach $\mathcal{A}$-bimodules.
 Assume that there exists a mapping $\varphi : M_1^n \rightarrow \mathbb {R}_+$
 for which a mapping $f : M_1 \rightarrow M_2$ satisfies the functional inequality
 \begin{eqnarray*}
&&\left\|n \sum_{1\le i<j\le n} f(ux_i-ux_j)
-\sum_{i=1}^{n} \hat{u} f\left(\sum_{j =1}^n x_j -n x_i\right)\right\| \\\nonumber
 &&\quad \le \varphi(x_1,\cdots,x_n), \quad \forall x_i \in M_1, \forall u \in \mathcal{A}(|u|=1),
\end{eqnarray*}
and the series (\ref{aqf3-1-1}) converges for all $x_i \in M_1$, $i=1,\cdots,n$.
If either $f$ is measurable or $f(tx)$ is continuous in $t \in \Bbb R$ for each
fixed $x \in M_1,$
then there exists a unique $\mathcal{A}_{sa}$-quadratic mapping $Q: M_1 \rightarrow M_2,$
defined by
$ Q(x)=\lim\limits_{m\rightarrow \infty} { f((n-1)^m x)  \over (n-1)^{2m}},$
which satisfies
the equation (\ref{fe3}) and the inequality (\ref{aqf3-1-2})
for all $x \in M_1.$
\end{thm}

%%%%%%%%%%%%%%%%%%%%%%%%%%%   Thm    %%%%%%%%%%%%%%%%%%%%%%%%%%%%%%%%%%%%%%%%%%%%%%%%%%

{\it Proof.}
By the same reasoning as the proof of Theorem \ref{aqf5}, it follows from $u=1 \in \mathcal{A}(|u|=1)$
 that there exists a unique
$\Bbb R$-quadratic mapping $Q : M_1 \rightarrow M_2,$
defined by
$ Q(x)=\lim\limits_{m \rightarrow \infty} { f((n-1)^m x)  \over (n-1)^{2m}},$
which satisfies the equation (\ref{fe3}) and the inequality (\ref{aqf3-1-2}).
By the similar manner to the proof of Theorem \ref{aqf5}
we obtain that $Q(u x)=\hat{u}Q(x)$ for all $x \in M_1$ and each $u \in \mathcal{A} (|u|=1).$
The last relation is also true for $u=0.$  Since $Q$ is $ \mathbb{R}$-quadratic,
for each element $a(a\ne 0) \in \mathcal{A}$
\begin{eqnarray*}
 Q(ax)  &=& Q\left( |a| \frac{a}{|a|}x \right)  = |a|^2 Q\left( \frac{a}{|a|}x \right)
 = |a|^2\frac{\hat{a}}{|a|^2} Q(x)\\
& =&\hat{a}Q(x)
\end{eqnarray*}
for all $x \in M_1$ and for all $a \in \mathcal{A}.$
 So the unique $\mathbb{R}$-quadratic mapping
$Q$ is also $\mathcal{A}_{sa}$-quadratic, as
desired. This completes the proof.
\hfill$\square$

%%%%%%%%%%%%%%%%%%%%%%%%%%%%%%%%%%%%%         Thm       %%%%%%%%%%%%%%%%%%%%%%%%

\begin{thm}\label{aqf9}
Let $M_1$ and $M_2$ be quasi-Banach $\mathcal{A}$-bimodules.
 Assume that there exists a mapping $\varphi : M_1^n \rightarrow \mathbb {R}_+$
 for which a mapping $f : M_1 \rightarrow M_2$ satisfies the functional inequality
 \begin{eqnarray*}
&&\left\|n \sum_{1\le i<j\le n} f(ux_i-ux_j)
-\sum_{i=1}^{n} \hat{u} f\left(\sum_{j =1}^n x_j -n x_i\right)\right\| \\\nonumber
 &&\quad \le \varphi(x_1,\cdots,x_n), \quad \forall x_i \in M_1, \forall u \in \mathcal{A}(|u|=1),
\end{eqnarray*}
and the series (\ref{aqf3-23}) converges for all $x_i \in M_1$, $i=1,\cdots,n$.
If either $f$ is measurable or $f(tx)$ is continuous in $t \in \Bbb R$ for each
fixed $x \in M_1,$
then there exists a unique $\mathcal{A}_{sa}$-quadratic mapping $Q: M_1 \rightarrow M_2,$
defined by
$Q(x)=\lim\limits_{m \rightarrow \infty}  (n-1)^{2m} f\Big( \frac{x}{(n-1)^m}  \Big),$
which satisfies
the equation (\ref{fe3}) and the inequality (\ref{aqf3-24})
for all $x \in M_1.$
\end{thm}

%%%%%%%%%%%%%%%%%%%%%%%%%%%%%%%%%====   Reference   =======%%%%%%%%%%%%%%%%%%%%%%%%%%%%

%\bibliographystyle{amsplain}

\end{document}